\documentclass[11pt]{article}
\usepackage[a4paper,margin=1in]{geometry}
\usepackage{amsmath,amssymb,amsthm,mathtools}
\usepackage[hidelinks]{hyperref}
\newtheorem{theorem}{Theorem}[section]

\newcommand{\PN}{\mathcal P_N}

\newcommand{\Q}{\mathbb Q}
\newcommand{\N}{\mathbb N}

\title{Large-prime residue bounds for reciprocal fractional-part sums}
\author{Ahmadreza Azimifard\thanks{Harmonic Research \& Technologies, LLC, New York, NY, USA.
  Email: \href{mailto:afard@harmonicrt.com}{afard@harmonicrt.com}}}
\hypersetup{pdftitle={Large-prime residue bounds for reciprocal fractional-part sums},
  pdfauthor={Ahmadreza Azimifard}}
\date{}

\begin{document}
\maketitle

\begin{abstract}
We record an explicit bound for the unrestricted integer-witness
maximum of $F_N(u)=\sum_{n\le N}\mu(n)\{u/n\}$.
The progression average for a reciprocal sawtooth and the Chinese
remainder theorem give a weighted finite bound. For M\"obius coefficients
the same bound follows from Caceres's primitive-sawtooth decomposition.
With $A_k=\sum_{m\le k}\mu(m)/m$, the prime number theorem gives
\[
 \liminf_{N\to\infty}\frac{\log N}{N}\sup_{u\ge0}|F_N(u)|
 \ \ge\ \frac12\sum_{k\ge1}\frac{|A_k|}{k(k+1)}
 \ \ge\ \frac{109}{360},
\]
where the supremum is over nonnegative integers. We give a direct
residue proof for rational coefficients and document its finite
formalization in Lean.
\end{abstract}

\section{Introduction}

For natural numbers $N,u$, consider
\[
 F_N(u)=\sum_{1\le n\le N}\mu(n)\{u/n\},
\]
where $\mu$ is the M\"obius function and $\{x\}$ is the fractional part.
We hold $N$ fixed while choosing $u$. The resulting integer witness may
range over a full common period of the denominators. The lower bound in
Theorem~\ref{thm:asymptotic} rules out a uniform $O(N^{1-\delta})$
upper bound for fixed $\delta>0$ with this choice of witnesses.
The prime number theorem enters only in the passage from the finite
inequality to the liminf \cite{MVI}.

\paragraph{Relation to earlier work.}
Averaging one sawtooth on a complete progression is the degree-one
case of Raabe's multiplication formula \cite{RaabeDLMF}.
Balakrishnan and P\'etermann \cite[Sec.~4, Eqs.~(4.9) and (4.11)]{BP1996}
use the same gcd-dependent averaging step for a weighted reciprocal
sawtooth sum with a moving cutoff. The complete residue cycles used
here give its exact finite version.
Montgomery used progressions on which a weighted sawtooth sum has a bias
\cite{Montgomery1987}. B\'ar\'any, Martin, Naslund and Robins
\cite[Section~3, especially Lemma~3.5]{BMNR2020} give this argument for
a totient error term whose corresponding sum has weights $\mu(d)/d$
and cutoff $d\le x$ at the evaluation point $x$.
Lemke Oliver and Soundararajan \cite[Proposition~3.1 and Section~5]{LOS}
isolate a prime occurring in just one denominator of a joint sawtooth
moment and study large values through periodic averages.

Caceres \cite[Sec.~10]{Caceres2026} decomposes the centered form of the same
M\"obius sum into primitive sawtooths. In Proposition~10.5, the
coefficient for a prime conductor $p>\sqrt N$ is
$-A_{\lfloor N/p\rfloor}$. The half-open sawtooth convention extends
the finite decomposition to integer arguments by right continuity.
Averaging on a complete integer progression with step equal to the
product of the primes in $(\sqrt N,N]$ keeps the constant and those
prime components. For a squarefree conductor $q$ with a prime
$\ell$ outside the step, the terms indexed by $t$ and $\ell t$ with
$\ell\nmid t$ cancel in the averaged divisor formula. The CRT selects
the endpoints in each surviving prime coordinate. Thus the M\"obius
bound is an explicit consequence of that decomposition. The continuous
mean-square identity in Caceres's Theorem~10.2 involves the same
GCD quadratic form studied in \cite{AzimifardGCD}.

Below we give a direct proof of the progression and range bounds for
arbitrary rational coefficients. The included Lean sources check the
finite part. We then sum over prime intervals to obtain the stated
liminf bound.

\section{Definitions and finite statements}

Throughout, $\N=\{0,1,2,\ldots\}$, sums over $1\le n\le N$ are empty
when $N=0$, and the coefficients $a_n$ are rational.  Set
\[
 \PN=\{p\le N: p\text{ prime},\ N<p^2\},\qquad
 Q_N=\prod_{p\in\PN}p,\qquad H=N!,
\]
and define, for $u,b\in\N$,
\begin{align*}
 F_N(a;u)&=\sum_{1\le n\le N}a_n\frac{u\bmod n}{n},\\
 C_p(a;N)&=\sum_{1\le m\le\lfloor N/p\rfloor}\frac{a_{mp}}m,\\
 \mathcal A_N(a;b)&=\frac1{N!}\sum_{0\le k<N!}F_N(a;b+kQ_N),\\
 \mathcal R_N(a)&=\sum_{p\in\PN}\left(1-\frac1p\right)|C_p(a;N)|.
\end{align*}
For positive $n$ and integer $u$, $(u\bmod n)/n=\{u/n\}$.

\begin{theorem}\label{thm:t0}
For every $a:\N\to\Q$ and $N,b,b'\in\N$,
\[
 \mathcal A_N(a;b)-\mathcal A_N(a;b')
 =\sum_{p\in\PN}\frac{(b\bmod p)-(b'\bmod p)}p\,C_p(a;N).
\]
\end{theorem}

\begin{theorem}\label{thm:t1}
(i) For every $a:\N\to\Q$ and $N\in\N$, there are
$b_+,b_-\in\N$ such that
\[
 \mathcal A_N(a;b_+)-\mathcal A_N(a;b_-)=\mathcal R_N(a).
\]
(ii) There are $u_+,u_-\in\N$ with
$\mathcal R_N(a)\le F_N(a;u_+)-F_N(a;u_-)$.
(iii) There is $u\in\N$ with
\[
 \frac12\mathcal R_N(a)\le |F_N(a;u)|.
\]
\end{theorem}

Put $A_k=\sum_{1\le m\le k}\mu(m)/m$.  The following statements use
$a_n=\mu(n)$, regarded as rational numbers.

\begin{theorem}\label{thm:t2}
(i) For $N\in\N$ and $p\in\PN$,
$C_p(\mu;N)=-A_{\lfloor N/p\rfloor}$.
(ii) For every $N\in\N$ there is
$u\in\N$ such that
\[
 \frac12\sum_{p\in\PN}\left(1-\frac1p\right)
 \left|A_{\lfloor N/p\rfloor}\right|\le |F_N(\mu;u)|.
\]
\end{theorem}

\begin{theorem}\label{thm:t3}
For every $N\in\N$ there is $u\in\N$ such that
\[
 \frac12\sum_{\substack{p\in\PN\\N<2p}}
 \left(1-\frac1p\right)\le |F_N(\mu;u)|.
\]
\end{theorem}

\section{Asymptotic statements}

Write $U_N=\sup_{u\in\N}|F_N(\mu;u)|$.

\begin{theorem}\label{thm:asymptotic}
(i) Using the prime number theorem,
\[
 U_N\ge \left(\frac14+o(1)\right)\frac{N}{\log N}.
\]
(ii) Moreover,
\[
 \liminf_{N\to\infty}\frac{\log N}{N}U_N
 \ge C_\mu:=\frac12\sum_{k\ge1}\frac{|A_k|}{k(k+1)}
 \ge\frac{109}{360}.
\]
These asymptotic assertions are paper-level statements and are not
part of the Lean formalization described in Section~\ref{sec:formal}.
\end{theorem}

\section{Proofs}

We give an elementary argument.  The discussion of
the partition uses $N\ge2$; the Lean declarations above also cover
$N=0,1$.

\begin{proof}[Proof of Theorems~\ref{thm:t0} and \ref{thm:t1}]
Two distinct primes in $\PN$ have product greater than $N$, so an
integer $n\le N$ is divisible by at most one of them.  If none divides
$n$, then $\gcd(Q_N,n)=1$.  If one does, write $n=mp$; then
$m<p$, and no prime of $\PN$ divides $m$.

When $d\mid H$ and $\gcd(c,d)=1$, the residues
$b+kc\pmod d$ for $0\le k<H$ visit each residue modulo $d$ exactly
$H/d$ times: each block of length $d$ is permuted by multiplication by
$c$.  For $n$ divisible by no prime of $\PN$, this gives average
$(n-1)/(2n)$ for $\{(b+kQ_N)/n\}$, independent of $b$.
For $n=mp$ with $p\in\PN$, put $r=b\bmod p$.  The residue modulo
$p$ stays fixed and the residue modulo $m$ is uniform.  The possible
residues modulo $mp$ are $r,r+p,\ldots,r+(m-1)p$, each equally often.
Their normalized average is
\[
 \frac{r}{mp}+\frac{m-1}{2m}.
\]
Grouping the finite sum by its unique large prime therefore expresses
$\mathcal A_N(a;b)$ as a term independent of $b$ plus
$\sum_{p\in\PN}(b\bmod p)C_p(a;N)/p$.  Taking a difference proves Theorem~\ref{thm:t0}.

The Chinese remainder theorem allows the residues modulo the distinct
primes in $\PN$ to be selected independently.  For $b_+$ choose residue
$p-1$ if $C_p\ge0$ and $0$ otherwise; for $b_-$ choose the opposite.
Theorem~\ref{thm:t0} then gives Theorem~\ref{thm:t1}(i).  An average is at most the
maximum of its terms and at least their minimum, giving witnesses
$u_+,u_-$ with the stated difference.  The triangle inequality yields
one witness with absolute value at least half that difference.
\end{proof}

\begin{proof}[Proof of Theorems~\ref{thm:t2} and \ref{thm:t3}]
For $p\in\PN$ and $m\le N/p$, one has $m<p$; hence $p\nmid m$ and
$\mu(mp)=-\mu(m)$.  This proves Theorem~\ref{thm:t2}(i). Theorem~\ref{thm:t1}(iii)
then gives Theorem~\ref{thm:t2}(ii).
For primes in the filtered sum of Theorem~\ref{thm:t3}, $p\le N<2p$ implies
$\lfloor N/p\rfloor=1$ and $A_1=1$.  Retaining those nonnegative terms
in Theorem~\ref{thm:t2}(ii) proves Theorem~\ref{thm:t3}.
\end{proof}

\begin{proof}[Proof of Theorem~\ref{thm:asymptotic}]
For $N\ge4$, every prime in $(N/2,N]$ belongs to $\PN$.  The prime
number theorem gives
$\pi(N)-\pi(N/2)\sim N/(2\log N)$, while the reciprocal-prime
correction in that interval is negligible.  Theorem~\ref{thm:t3} gives Theorem~\ref{thm:asymptotic}(i).

Fix $K$.  For $1\le k\le K$, the primes with
$\lfloor N/p\rfloor=k$ lie in $(N/(k+1),N/k]$; when
$N>(K+1)^2$ all of them belong to $\PN$.  The prime number theorem gives
\[
 \pi(N/k)-\pi(N/(k+1))
 \sim \frac{N}{k(k+1)\log N},
\]
and the reciprocal-prime correction is $O_k(1)$.  Theorem~\ref{thm:t2}(ii) then yields
\[
 \liminf_{N\to\infty}\frac{\log N}{N}U_N
 \ge \frac12\sum_{k=1}^{K}\frac{|A_k|}{k(k+1)}.
\]
The terms are nonnegative, and $|A_k|\le H_k\ll\log k$ for large $k$,
so the series converges and taking all $K$ proves the liminf bound.
Finally, the exact values
$A_1=1$, $A_2=1/2$, $A_3=A_4=1/6$ yield
\[
 \frac12\left(\frac12+\frac{1/2}{6}
                    +\frac{1/6}{12}+\frac{1/6}{20}\right)
 =\frac{109}{360},
\]
which gives the numerical lower bound in Theorem~\ref{thm:asymptotic}(ii).
\end{proof}

\section{Formalization}\label{sec:formal}

The finite statements in Theorem~\ref{thm:t0} and Theorems~\ref{thm:t1},
\ref{thm:t2} and \ref{thm:t3} were checked in Lean 4.28.0 with Mathlib
4.28.0, as pinned in the ancillary files. The formalization uses
$\Q$-valued coefficients, residues of natural-number witnesses and
factorial averaging. It does not formalize a real supremum or the
asymptotic statements. The 12 checked declarations use only
\texttt{propext}, \texttt{Classical.choice} and \texttt{Quot.sound}.
See \texttt{anc/README.txt} for the correspondence. The ancillary sources
and project configuration were not rebuilt for this submission.

\paragraph{Acknowledgment.}
Automated research and drafting tools were used to derive and prepare
this note. The scope of the Lean checks is specified above.

{\small
\bibliographystyle{plain}
\bibliography{refs}

@book{MVI,
  author = {Montgomery, H. L. and Vaughan, R. C.},
  title = {Multiplicative Number Theory {I}: Classical Theory},
  series = {Cambridge Studies in Advanced Mathematics}, volume = {97},
  publisher = {Cambridge University Press}, year = {2007}
}

@article{BMNR2020,
  author = {B{\'a}r{\'a}ny, Imre and Martin, Greg and Naslund, Eric and Robins, Sinai},
  title = {Primitive points in rational polygons},
  journal = {Canadian Mathematical Bulletin}, volume = {63}, number = {4},
  year = {2020}, pages = {850--870}
}

@misc{Caceres2026,
  author = {Caceres, Pedro},
  title = {{The Riemann hypothesis as a margin: a two-sided criterion verified to
           $2\cdot 10^{29}$, the fourth-moment anatomy of M\"obius cancellation,
           and one constant seen from every side}},
  year = {2026},
  note = {Hexagon preprint 2610.00688v3}
}

@misc{RaabeDLMF,
  author = {{NIST Digital Library of Mathematical Functions}},
  title = {Multiplication formulas for {Bernoulli} polynomials},
  note = {Section 24.4(v), Eq.~24.4.18}
}

@article{Montgomery1987,
  author = {Montgomery, Hugh L.},
  title = {Fluctuations in the mean of {Euler}'s phi function},
  journal = {Proceedings of the Indian Academy of Sciences, Mathematical Sciences},
  volume = {97}, year = {1987}, pages = {239--245}
}

@article{LOS,
  author = {Lemke Oliver, Robert J. and Soundararajan, Kannan},
  title = {The distribution of consecutive prime biases and sums of sawtooth random variables},
  journal = {Mathematical Proceedings of the Cambridge Philosophical Society},
  volume = {168}, number = {1}, pages = {149--169}, year = {2020}
}

@misc{AzimifardGCD,
  author = {Azimifard, Ahmadreza},
  title = {A power-saving asymptotic for a {M\"obius}-weighted {GCD} quadratic form},
  year = {2026}, note = {arXiv:2611.00525}
}

@article{BP1996,
  author={Balakrishnan, Ukath-Varyat and P{\'e}termann, Yves-Fran{\c{c}}ois},
  title={The {Dirichlet} series of {$\zeta(s)\zeta^\alpha(s+1)f(s+1)$}:
         On an error term associated with its coefficients},
  journal={Acta Arithmetica}, volume={75}, number={1}, year={1996}, pages={39--69}
}
}
\end{document}